\documentclass[11pt]{amsart}
\usepackage[utf8]{inputenc}
\usepackage[total={6.5in,9in},top=1in, left=1in, right=1in, bottom=1in]{geometry}
\usepackage{fancyhdr}
\usepackage{color,amsmath,amssymb,mathtools,graphicx}
\usepackage{amsfonts}
\usepackage{libertine}
\usepackage[libertine]{newtxmath}
\usepackage{bbold}
\usepackage{epsfig,fancyhdr}
\usepackage{dsfont} 
\usepackage{graphicx,float}
\usepackage{graphicx}
\usepackage{epsfig,color,fancyhdr,setspace}
\usepackage{epstopdf}
\usepackage{wrapfig}
\usepackage{soul}
\usepackage{import}
\usepackage{lineno}
\usepackage{stackrel}

\usepackage[hidelinks,pagebackref,pdftex]{hyperref}

\usepackage[abs]{overpic}
\usepackage[center]{caption}
\usepackage{subcaption}
\usepackage{tikz-cd}
\usepackage{enumitem}
\counterwithin{figure}{section}
\usepackage{tikz-cd}
\usepackage{wasysym}
\usepackage{marvosym}
\usepackage{extarrows}
\usepackage{amsmath}
\usepackage{graphicx,amsmath,mathtools,float}
\usepackage{setspace}
\graphicspath{{./images/}}
\usepackage{stackrel}
\usepackage{yfonts}
\usepackage{enumitem, hyperref}\makeatletter
\def\namedlabel#1#2{\begingroup
    #2%
    \def\@currentlabel{#2}%
    \phantomsection\label{#1}\endgroup
}
\makeatother
\usepackage{type1cm} 
\usepackage{blindtext}
\usepackage{scalerel,stackengine}

\usepackage{tikz}

\usepackage{mathtools}
\usepackage{amssymb}
\usepackage{tikz}
\usetikzlibrary{shapes,snakes}

\newcommand{\KP}[1]{%
	\begin{tikzpicture}[baseline=-\dimexpr\fontdimen22\textfont2\relax]
		#1
	\end{tikzpicture}%
}

\newcommand{\KPB}{%
	\KP{
		\draw[color=gray,thick] (-0.3,0.3) -- (0.3,-0.3);
		\draw[color=gray,thick] (-0.3,-0.3) -- (-0.05,-0.05);
		\draw[color=gray,thick] (0.05,0.05) -- (0.3,0.3);
	}%
}
\newcommand{\KPC}{%
	\KP{%
		\draw[color=gray,thick] (-0.3,0.3) .. controls (0,-0.05) .. (0.3,0.3);
		\draw[color=gray,thick] (-0.3,-0.3) .. controls (0,0.05) .. (0.3,-0.3);
	}%
}
\newcommand{\KPD}{%
	\KP{%
		\draw[color=gray,thick] (-0.3,-0.3) .. controls (0.05,0) .. (-0.3,0.3);
		\draw[color=gray,thick] (0.3,-0.3) .. controls (-0.05,0) .. (0.3,0.3);
	}%
}

\newcommand\quotient[2]{
	\mathchoice
	{
		\text{\raise1ex\hbox{$#1$}\Big/\lower1ex\hbox{$#2$}}%
	}
	{
		#1\,/\,#2
	}
	{
		#1\,/\,#2
	}
	{
		#1\,/\,#2
	}
}

\usepackage{tocbasic}
\newtheorem{theorem}{Theorem}[section]

\newtheorem{definition}[theorem]{Definition}

\newtheorem*{conjecture*}{Alternation Conjecture}
\newtheorem*{conj*}{Generalized Kauffman-Harary Conjecture}
\newtheorem*{conje*}{Alternate Forms of the Generalized Kauffman-Harary Conjecture}

\title{A Glimpse of the Khovanov Homology of $\boldsymbol{T(2,n)}$ Via Long Exact Sequence}

\author{Gabriel Montoya-Vega}
\address{Department of Mathematics, University of Puerto Rico at R\'io Piedras, San Juan, PR, USA }
\email{\textcolor{blue}{gabrielmontoyavega@gmail.com $|$ gabriel.montoya@upr.edu}}



\begin{document}

\begin{abstract}
	Khovanov homology is a powerful link invariant: a categorification of the Jones polynomial that enjoys a rich and beautiful algebraic structure. This homology theory has been extensively studied and it has become an ubiquitous topic in contemporary knot theory research. In the same spirit, the Kauffman skein relation, which allows to define the Kauffman bracket polynomial up to normalization of the unknot, can be categorified by means of a long exact sequence. In an expository style, in this article we present how to build Khovanov homology from the Kauffman bracket polynomial and construct its long exact sequence. Furthermore, we present a deviceful and practical way in which this long exact sequence can be used for the computation of the Khovanov homology of torus links of the type $T(2,n)$.
	
	\ 
	
	This article serves as a partial translation of a Spanish paper to be published on occasion of the \textit{Encuentro Internacional de Matem\'aticas} (International Meeting of Mathematics) celebrated at the Universidad del Atl\'antico in Barranquilla, Colombia in November 2023. This paper offers a first look into the world of Khovanov homology by constructing it from the Kauffman bracket polynomial, as it was first done by Oleg Viro. Moreover, it gives the reader references for further studies from leading experts such as D. Bar-Natan, M. Khovanov, S. Mukherjee, J. Przytycki, and A. Shumakovitch, among others. In particular, one of the main objectives in publishing this article (and this partial translation) is to popularize research in knot theory, more specifically on Khovanov homology in Colombia, and Latin-America in general, acting as a language bridge given that most of the literature is in English.

\end{abstract}

\maketitle

\tableofcontents

\section{Introduction}

Khovanov homology (KH) offers a nontrivial generalization of the Jones polynomial and the Kauffman bracket polynomial of links in $\mathbb{R}^{3}$. A more powerful invariant than the Jones polynomial, it has been extensively studied over the last two decades. The essence of KH is that a bigraded chain complex is associated to a link, in such a way that the homology of the complex is a link invariant. Furthermore, the graded Euler characteristic of the chain complex is the Jones polynomial.

\ 

This article is organized as follows. In Section \ref{KHconstruction} we present how to build Khovanov homology from the Kauffman bracket polynomial and show how this approach connects with the original construction. In Section \ref{LESconstruction}, the long exact sequence of KH is constructed and in Section \ref{KHtoruslinks} we use it to compute the KH of torus links $T(2,n)$. Section \ref{examplesKHtables} consists of two examples illustrating the KH of two specific knots. Finally, in Section \ref{Futuredirec} we give some possibilities for further research paths in the context of different families of links.

\section{Khovanov homology via Kauffman bracket polynomial}	\label{KHconstruction}

In this section we construct KH following Oleg Viro's approach \cite{Vir1, Vir2}. Recall that the Kauffman bracket polynomial (KBP) gives a simple approach to the Jones polynomial \cite{Kau}. Hence, it is natural to first categorify the KBP and then connect this theory to Khovanov's original \cite{Kho1}. The construction of KH that we present here is also referred to as the unoriented framed version of KH.

\begin{definition}\ \label{KBPdefi}
	\begin{itemize}
		\item [(i)] The (unreduced) \textbf{Kauffman bracket polynomial} is a function from the set of unoriented link diagrams $\mathcal{D}$ to Laurent polynomials with integer coefficients in the variable $A$, $\left[ \ \right]: \mathcal{D} \longrightarrow \ \mathbb{Z}[A^{\pm1}]$. The polynomial is characterized by the rules $[\emptyset]=1$, $[\bigcirc]=(-A^{2}-A^{-2})$, and the skein relation:
		$$\left[ \KPB \right]=
		A\left[ \KPC \right] + A^{-1} \left[ \KPD \right].$$
		\item [(ii)]	Let $D$ be an unoriented link diagram and let $cr(D)$ be its crossings set. A \textbf{Kauffman state} s, of $D$, is a function $s: cr(D)\longrightarrow \{A, B\}$. This function is understood as an assignment of a \textbf{marker} to each crossing according to the convention illustrated in Figure \ref{fig:kaufstatefunction}. Denote by \textit{KS} the set of all Kauffman states. Every marker yields a natural \textbf{smoothing} of the crossing as shown in  Figure \ref{fig:kaufstatefunction}. 
		\begin{figure}[ht]
			\centering
			\includegraphics[width=0.95\linewidth]{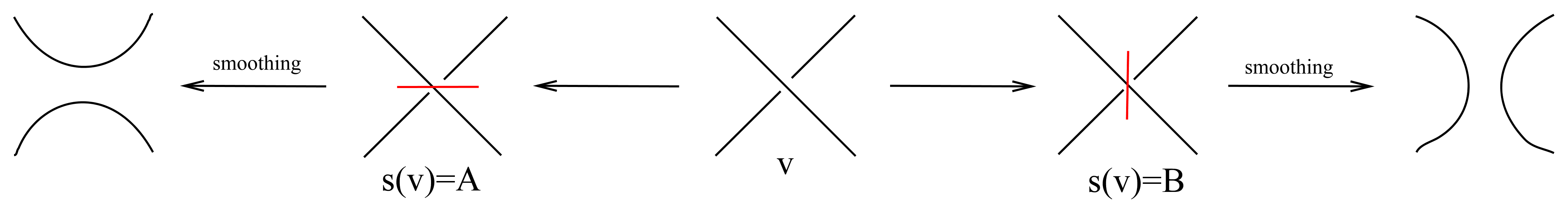}
			\caption{Markers at a crossing v of $D$ and their corresponding smoothing.}
			\label{fig:kaufstatefunction}
		\end{figure}
	\end{itemize}
	
\end{definition}

The KBP of a link diagram $D$ is given by the state sum  formula: 
$$[D]=\sum_{s \ \in \ KS}^{}A^{\mid s^{-1}(A) \mid - \mid s^{-1}(B) \mid }(-A^{2}-A^{-2})^{\mid D_{s} \mid},$$
where $D_{s}$ denotes the system of circles obtained after smoothing all crossings of $D$ according to the markers of $s$, and $|D_{s}|$ denotes the number of circles in the system.

\
 
By having the state sum formula from the unreduced version of the Kauffman bracket polynomial, it is possible to have an one-to-one correspondence between the circles in $D_{s}$ and the factors $(-A^{2}-A^{-2})$. This leads to the idea of an enhanced Kauffman state, as defined below.	

\begin{definition}\label{EnhancedStatesDefi}
	An \textbf{enhanced Kauffman state} $S$ of $D$ is a Kauffman state $s$ together with a function  $\varepsilon: D_{s}  \longrightarrow \{+, -\}$, assigning to each circle of $D_{s}$ a positive or a negative sign.
\end{definition}	

Denote by $EKS$ the set of all enhanced Kauffman states. Notice that for a Kauffman state $s$ there are $2^{\mid D_{s} \mid}$ enhanced Kauffman states. Then the KBP can be written as a sum of monomial terms coming from the enhanced Kauffman states as follows:
$$[D]=\sum_{S \ \in \ EKS}^{}(-1)^{\mid D_{s} \mid}A^{\mid s^{-1}(A) \mid - \mid s^{-1}(B) \mid }(A^{2})^{\mid \varepsilon^{-1}(+) \mid - \mid \varepsilon^{-1}(-) \mid }.$$ 
Letting $\sigma (s)=|s^{-1}(A)|-|s^{-1}(B)|$ and $\tau (S)=|\varepsilon^{-1}(+)|-|\varepsilon^{-1}(-)|$, we obtain the following formula:
$$[D]=\sum_{S \ \in \ EKS}^{}(-1)^{\mid D_{s} \mid}A^{\sigma (s)+2\tau (S)},$$	
which is the enhanced Kauffman state sum formula for the unreduced KBP.

\

The enhanced Kauffman states form a basis for the chain groups of the Khovanov chain complex. We now define the bigrading on EKS, the chain groups, and the boundary maps.

\begin{figure}[h]
	\centering
	\includegraphics[width=0.55\linewidth]{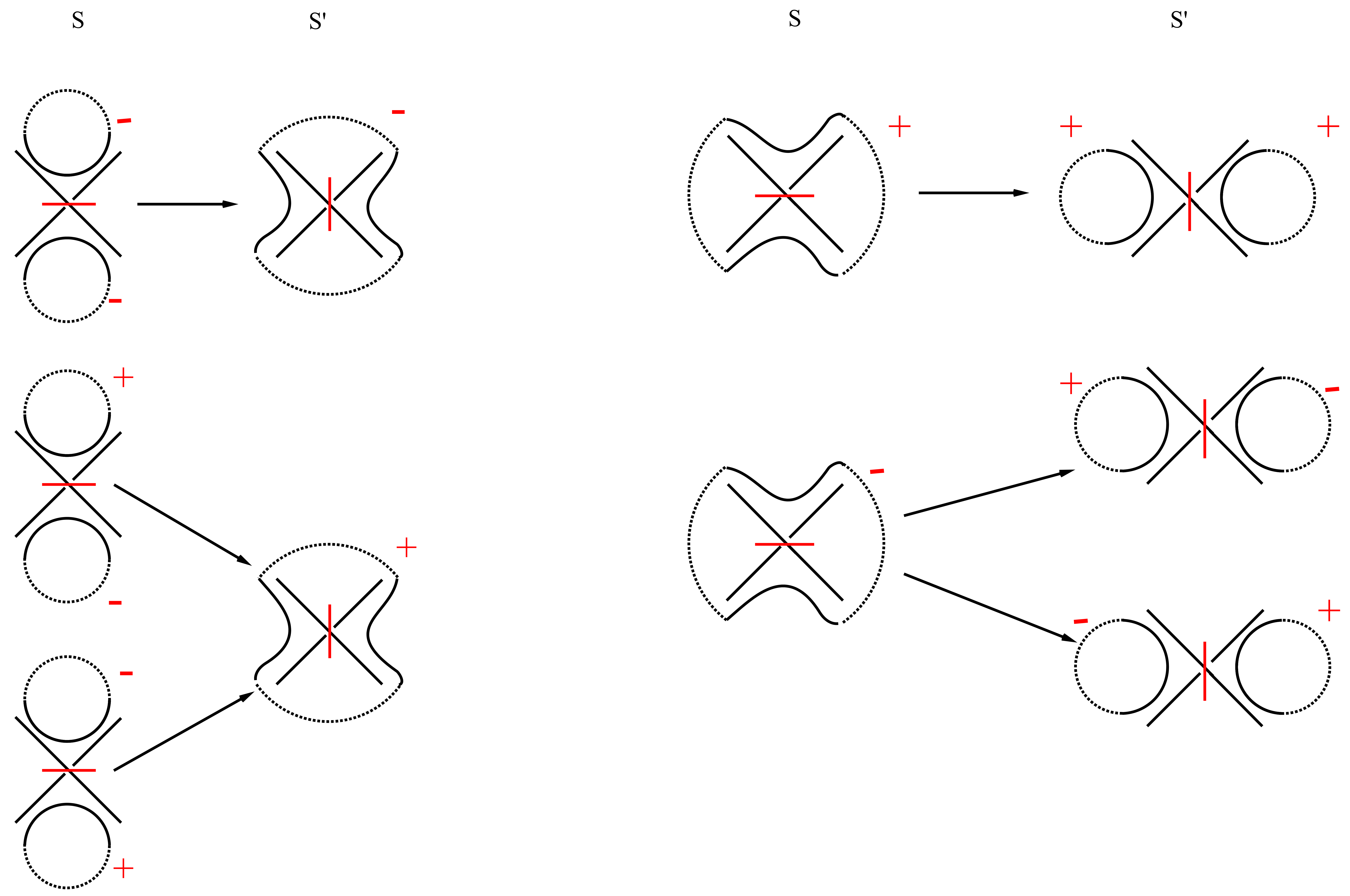}
	\caption{Signs of circles after fusion (left) and split (right) of components.}
	\label{fig:cases-merged}
\end{figure}

\begin{definition}\ \label{basis}
	
	\begin{enumerate}
		\item [(i)]The \textbf{bidegree} on the enhanced Kauffman states is defined as the following set: $$\mathcal{S}_{a,b}(D)=\mathcal{S}_{a,b}=\left\lbrace  S \ \in \ EKS \ \mid \  a=\sigma (s), \  b=\sigma (s)+2\tau (S) \right\rbrace. $$
		\item [(ii)] The \textbf{chain groups} $\mathcal{C}_{a,b}(D)=\mathcal{C}_{a,b}$, are defined to be the free abelian groups with basis $\mathcal{S}_{a,b}(D)=\mathcal{S}_{a,b}$, i.e. $\mathcal{C}_{a,b}=\mathbb{Z}\mathcal{S}_{a,b}$. Therefore, $\mathcal{C}(D)=\displaystyle\bigoplus_{a,b \ \in \ \mathbb{Z}}\mathcal{C}_{a,b}(D)$ is a bigraded free abelian group.
		\item [(iii)] For a link diagram $D$ we define the \textbf{chain complex} $\mathscr{C}(D)= \left\lbrace \left( C_{a,b}, \partial_{a,b} \right) \right\rbrace  $,
		where 
		the \textbf{differential map} $\partial_{a,b}: \mathcal{C}_{a,b} \longrightarrow \mathcal{C}_{a-2,b}$ is defined by
		$$\partial_{a,b}(S)= \sum_{S' \ \in \ \mathcal{S}_{a-2,b} }^{}(-1)^{t(S,S')}(S,S')S'.$$
	\end{enumerate}
	
\end{definition}
In the definition of the differential map above, $S \in \mathcal{S}_{a,b}$ and $(S,S')$ is the incidence number of $S$ and $S'$. Specifically, either $(S,S')=0$ or $(S, S')=1$, and it is equal to $1$ if and only if the following conditions hold:
\begin{enumerate}
	\item $S$ and $S'$ are identical except at only one crossing, say $v$. Moreover $s(v)=A$ and $s'(v)=B$, where $S$ and $S'$ are enhancements of $s$ and $s'$, respectively. 
	\item $\tau (S')=\tau (S)+1$ and every component of $D_{s}$ not interacting with the crossing $v$ keeps its sign for $D_{s'}$. 
\end{enumerate}	
Condition (1) reflects the fact that the value of $\sigma $ is decreasing by 2, while condition (2) indicates that either the number of negative signs decreases or the number of positive signs increases. Finally, $(-1)^{t(S,S')}$ requires an ordering of the crossings in the link diagram $D$. $t(S,S')$ is defined as the number of crossings of $D$ with $B$ markers in $S$ bigger than $v$ in the chosen ordering. This condition is sufficient for the differential to satisfy $\partial_{a-2,b}\circ \partial_{a,b}=0$. It is important to remark that the homology does not depend on the ordering of crossings \cite{Kho1}. Figure \ref{fig:cases-merged} illustrates the cases when the enhanced states $S$ and $S'$  are incident.

\begin{definition}\label{khdefi}
	The \textbf{Khovanov homology} of the diagram $D$ is defined to be the homology of the chain complex $\mathscr{C}(D)$:		$$H_{a,b}(D)=\dfrac{ker (\partial_{a,b}) }{im(\partial_{a+2,b})}.$$
\end{definition}

Certainly, the most important property of the framed version of Khovanov homology is its invariance under second and third Reidemeister moves. The following theorem summarizes the results. See \cite{Vir1, Vir2} for a proof.

\begin{theorem}\label{khinvariance}
	Let $D$ be a link diagram. The homology groups
	$$H_{a,b}(D)=\dfrac{ker (\partial_{a,b}) }{im(\partial_{a+2,b})},$$
	are invariant under Reidemeister moves of second and third type. Therefore, they are invariants of unoriented framed links. Moreover, the effect of the first Reidemeister move (positive or negative) $R_{1}$, is the shift in the homology, $H_{a,b}(R_{1+}(D))=H_{a+1,b+3}(D)$ and $H_{a,b}(R_{1-}(D))=H_{a-1,b-3}(D)$. These groups categorify the unreduced Kauffman bracket polynomial and are called the framed Khovanov homology groups.   
\end{theorem}

\subsection*{Classical Khovanov (co)homology }\label{classickh}

Khovanov originally associated a bigraded chain complex to an oriented link diagram whose homology is a link invariant. Let $\vec{D}$ be an oriented link diagram obtained after choosing an orientation on the unoriented link diagram $D$. Let $w(\vec{D})=w$ be its writhe. Then, the classical Khovanov (co)homology, $\mathcal{H}^{i,j}(\vec{D})$, and the framed version of KH, $H_{a,b}(D)$, are related by the following equalities: 
$$\mathcal{H}^{i,j}(\vec{D})=H_{w-2i,3w-2j}(D)=H_{a,b}(D)=\mathcal{H}^{\frac{w-a}{2},\frac{3w-b}{2}}(\vec{D}).  $$

\section{Long exact sequence of Khovanov homology}\label{LESconstruction}
In this section, we examine the fact that the skein relation used to define the KBP can be categorified by means of a long exact sequence \cite{Vir1, Vir2}. Then this long exact sequence is used to calculate the KH of torus links of type $T(2,n)$.

\subsection*{Construction}
Recall from Definition \ref{basis} that the set of enhanced Kauffman states $\mathcal{S}_{a,b}$ gives the basis for the free abelian groups  $\mathcal{C}_{a,b}=\mathbb{Z}\mathcal{S}_{a,b}$. Let $v$ be a fixed crossing of the link diagram $D$. Consider the sets 
$\mathcal{S}_{a,b}^{A,v}$ and $\mathcal{S}_{a,b}^{B,v}$ defined as follows:
$$\mathcal{S}_{a,b}^{A,v}=\left\lbrace S \ \in \ S_{a,b} \ \mid \ s(v)=A \right\rbrace \ \ \ \ \ and  \ \ \ \  \ \mathcal{S}_{a,b}^{B,v}=\left\lbrace S \ \in \ S_{a,b} \ \mid \ s(v)=B \right\rbrace,$$
this is, $\mathcal{S}_{a,b}^{A,v}$ consists of all enhanced states $S$ with bigrading $(a,b)$ having an $A$ marker at the crossing $v$; analogously, $\mathcal{S}_{a,b}^{B,v}$ consists of all enhanced states $S$ with bigrading $(a,b)$ having a $B$ marker at the crossing $v$. Then, it can be seen that 	
$$\mathcal{S}_{a,b} \ \ = \ \   \mathcal{S}_{a,b}^{A,v} \ \ \sqcup \ \ \mathcal{S}_{a,b}^{B,v}.$$

Continuing with the notation from the previous section, denote the free abelian groups generated by these sets as:
$$\mathbb{Z}\mathcal{S}_{a,b}^{A,v} \ \ = \ \ \mathcal{C}_{a,b}^{\includegraphics[width=0.04\textheight]{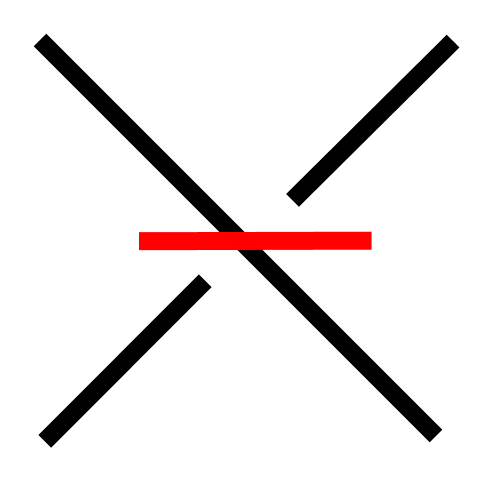}} \ \ \ \ \  and \ \ \ \ \ \mathbb{Z}\mathcal{S}_{a,b}^{B,v} \ \ = \ \ \mathcal{C}_{a,b}^{\includegraphics[width=0.04\textheight]{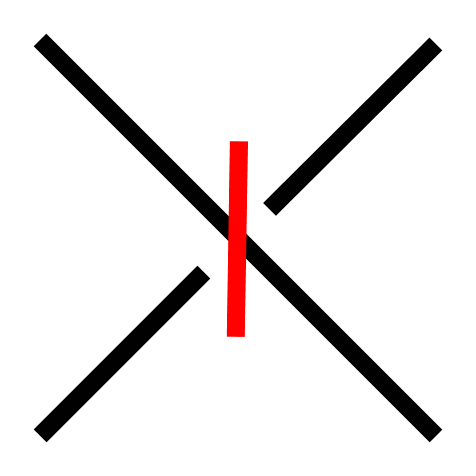}},$$
where $v= \textbf{\KPB}$. Consequently, at the level of groups we have that:
$$\mathcal{C}_{a,b} \ \ = \ \ \mathcal{C}_{a,b}^{\includegraphics[width=0.04\textheight]{positive_marker}} \ \ \oplus \ \ \mathcal{C}_{a,b}^{\includegraphics[width=0.04\textheight]{negative_marker}} \ \ \ or \ \  equivalently \ \ \ \mathbb{Z}\mathcal{S}_{a,b} \ \ = \ \ \mathbb{Z}\mathcal{S}_{a,b}^{A,v} \ \ \sqcup \ \ \mathbb{Z}\mathcal{S}_{a,b}^{B,v}.$$ 

Notice that the complex $( \mathcal{C}_{a,b}^{\includegraphics[width=0.04\textheight]{negative_marker}}, \partial_{a,b}) $ is a chain subcomplex of $\left( \mathcal{C}_{a,b}, \partial_{a,b}  \right) $, this is $\partial(\mathcal{C}_{a,b}^{\includegraphics[width=0.04\textheight]{negative_marker}}) \subset \mathcal{C}_{a-2,b}^{\includegraphics[width=0.04\textheight]{negative_marker}}$. In contrast, this is not necessarily true for $(\mathcal{C}_{a,b}^{\includegraphics[width=0.04\textheight]{positive_marker}}, \partial_{a,b})$ since $\partial_{a,b}$ may change the marker at the crossing from an $A$ marker to a $B$ marker, as discussed while constructing KH in Section \ref{KHconstruction}. The following short exact sequence of chain complexes can be written: 
$$0 \xrightarrow{\makebox[1.3cm]{}}\mathcal{C}_{a,b}^{\includegraphics[width=0.04\textheight]{negative_marker}} \xrightarrow{\makebox[1.3cm]{\textcolor{blue}{$\alpha$}}} \mathcal{C}_{a,b}^{\includegraphics[width=0.04\textheight]{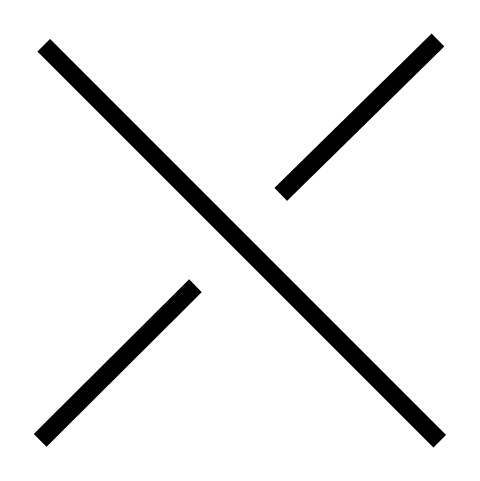}} \xrightarrow{\makebox[1.3cm]{\textcolor{blue}{$\beta$}}}  \quotient{\textstyle \mathcal{C}_{a,b}^{\includegraphics[width=0.04\textheight]{crossing}}}{\mathcal{C}_{a,b}^{\includegraphics[width=0.04\textheight]{negative_marker}}}  \xrightarrow{\makebox[1.3cm]{}} 0,$$
\noindent where the map $\alpha$ sends an enhanced Kauffman state of $D$ having the crossing $v$ with a $B$ marker to the enhanced Kauffman state of $D$ assigning a $B$ label to $v$, while the other crossings keep the markers and the signs of the circles are preserved. Similarly, the map $\beta$ sends an enhanced Kauffman state having a $B$ label at $v$ to zero, and sends each enhanced Kauffman state with $A$ label at $v$ to the enhanced Kauffman state of $D$ with the crossing $v$ given an $A$ label, while the other markers of crossings and signs of circles are preserved.

\

Observe that as a group $\quotient{\textstyle \mathcal{C}_{a,b}^{\includegraphics[width=0.04\textheight]{crossing}}}{\mathcal{C}_{a,b}^{\includegraphics[width=0.04\textheight]{negative_marker}}} \ \ = \ \ \mathcal{C}_{a,b}^{\includegraphics[width=0.04\textheight]{positive_marker}}$, and there is a chain complex connecting map coming from the construction of KH: 
$$\quotient{\textstyle \mathcal{C}_{a,b}^{\includegraphics[width=0.04\textheight]{crossing}}}{\mathcal{C}_{a,b}^{\includegraphics[width=0.04\textheight]{negative_marker}}} \xrightarrow{\makebox[3cm]{\textcolor{blue}{$\partial^{Conn}_{a,b}$}}} \mathcal{C}_{a-2,b}^{\includegraphics[width=0.04\textheight]{negative_marker}}.$$
As it is standard, this leads to a long exact sequence of homology:
\begin{equation*}
	\begin{split}
		\cdots \xrightarrow{\makebox[0.8cm]{}} H_{a,b}( \includegraphics[width=0.04\textheight]{negative_marker}) &   \xrightarrow{\makebox[1.3cm]{\textcolor{blue}{$\alpha_{\ast}$}}} H_{a,b}(\includegraphics[width=0.04\textheight]{crossing})\xrightarrow{\makebox[1.6cm]{\textcolor{blue}{$\beta_{\ast}$}}}H_{a,b}\left(\quotient{\textstyle \includegraphics[width=0.04\textheight]{crossing}} {\includegraphics[width=0.04\textheight]{negative_marker}}\right)  \\
		\xrightarrow{\makebox[1.8cm]{$  \textcolor{blue}{\left( \partial^{Conn}_{a,b} \right)_{\ast}}$}} H_{a-2,b}(\includegraphics[width=0.04\textheight]{negative_marker}) & \xrightarrow{\makebox[1.3cm]{\textcolor{blue}{$\alpha_{\ast}$}}}H_{a-2,b}(\includegraphics[width=0.04\textheight]{crossing}) \xrightarrow{\makebox[1.2cm]{\textcolor{blue}{$\beta_{\ast}$}}}H_{a-2,b}\left(\quotient{\textstyle \includegraphics[width=0.04\textheight]{crossing}} {\includegraphics[width=0.04\textheight]{negative_marker}}\right)  \\
		\xrightarrow{\makebox[1.8cm]{$ \textcolor{blue}{\left( \partial^{Conn}_{a,b} \right)_{\ast}}$}} H_{a-4,b}(\includegraphics[width=0.04\textheight]{negative_marker}) &  \xrightarrow{\makebox[1.3cm]{\textcolor{blue}{$\alpha_{\ast}$}}}H_{a-4,b}(\includegraphics[width=0.04\textheight]{crossing})\xrightarrow{\makebox[1.6cm]{\textcolor{blue}{$\beta_{\ast}$}}}\cdots.  \\ 
	\end{split}
\end{equation*}	

Observe that after smoothing the crossing, we have the following chain complex equalities:
$$\mathcal{C}_{a,b}^{\includegraphics[width=0.04\textheight]{negative_marker}} \ \ = \ \ \mathcal{C}_{a+1,b+1}^{\includegraphics[width=0.04\textheight]{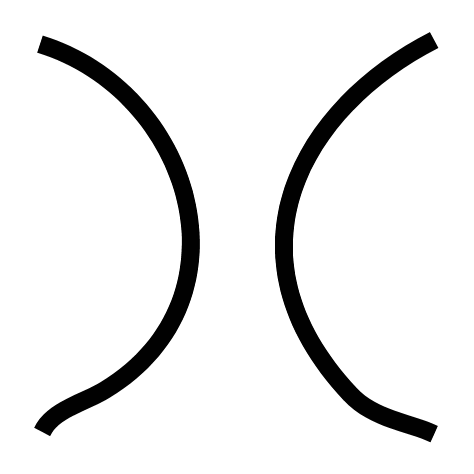}} \ \ \ \ \ \ and \ \ \ \ \ \ \quotient{\textstyle \mathcal{C}_{a,b}^{\includegraphics[width=0.04\textheight]{crossing}}}{\mathcal{C}_{a,b}^{\includegraphics[width=0.04\textheight]{negative_marker}}} \ \ = \ \ \mathcal{C}_{a-1,b-1}^{\includegraphics[width=0.04\textheight]{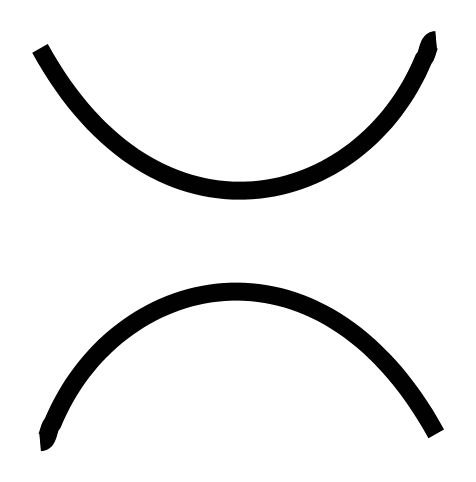}}. $$
Hence, the following short exact sequence of chain complexes of diagrams is obtained:
$$0 \xrightarrow{\makebox[1.5cm]{}}\mathcal{C}_{a+1,b+1}^{\includegraphics[width=0.04\textheight]{b_smoothing}} \xrightarrow{\makebox[1.5cm]{\textcolor{blue}{$\alpha$}}} \mathcal{C}_{a,b}^{\includegraphics[width=0.04\textheight]{crossing}} \xrightarrow{\makebox[1.5cm]{\textcolor{blue}{$\beta$}}}\mathcal{C}_{a-1,b-1}^{\includegraphics[width=0.04\textheight]{a_smoothing}}    \xrightarrow{\makebox[1.5cm]{}} 0,$$

\noindent which in turn, yields the following long exact sequence, known as the Long Exact Sequence of Khovanov homology:

\begin{equation}\label{LES}
	\begin{split}
		\cdots \xrightarrow{\makebox[1cm]{}} H_{a+1,b+1}(\includegraphics[width=0.035\textheight]{"b_smoothing"}) & \xrightarrow{\makebox[1cm]{\textcolor{blue}{$\alpha_{\ast}$}}} H_{a,b}(\includegraphics[width=0.035\textheight]{"crossing"})\xrightarrow{\makebox[1.3cm]{\textcolor{blue}{$\beta_{\ast}$}}}H_{a-1,b-1}(\includegraphics[width=0.035\textheight]{a_smoothing}) \\
		\xrightarrow{\makebox[2cm]{$  \textcolor{blue}{\left( \partial^{Conn}_{a,b} \right)_{\ast}}$}} H_{a-1,b+1}(\includegraphics[width=0.035\textheight]{"b_smoothing"}) & \xrightarrow{\makebox[1cm]{\textcolor{blue}{$\alpha_{\ast}$}}} H_{a-2,b}(\includegraphics[width=0.035\textheight]{"crossing"})\xrightarrow{\makebox[1cm]{\textcolor{blue}{$\beta_{\ast}$}}}H_{a-3,b-1}(\includegraphics[width=0.035\textheight]{"a_smoothing"}) \\
		\xrightarrow{\makebox[2cm]{$  \textcolor{blue}{\left( \partial^{Conn}_{a,b} \right)_{\ast}}$}} H_{a-3,b+1}(\includegraphics[width=0.035\textheight]{"b_smoothing"}) & \xrightarrow{\makebox[1cm]{\textcolor{blue}{$\alpha_{\ast}$}}}H_{a-4,b}(\includegraphics[width=0.035\textheight]{"crossing"})\xrightarrow{\makebox[1.3cm]{\textcolor{blue}{$\beta_{\ast}$}}} \cdots. \\
	\end{split}
\end{equation}

The following results follows directly from the previous construction.
\begin{theorem}\    \label{Coro1}
	\begin{enumerate}
		\item [(1)] If $H_{a+1,b+1}(\includegraphics[width=0.025\textheight]{"b_smoothing"})=0$ then $\beta_*: H_{a+1,b+1}(\includegraphics[width=0.025\textheight]{"crossing"}) \to H_{a,b}(\includegraphics[width=0.025\textheight]{"a_smoothing"})$ is a monomorphism.
		\item[(2)] If $H_{a-1,b+1}(\includegraphics[width=0.025\textheight]{"b_smoothing"})=0$ then $\beta_*: H_{a+1,b+1}(\includegraphics[width=0.025\textheight]{"crossing"}) \to H_{a,b}(\includegraphics[width=0.025\textheight]{"a_smoothing"})$ is an epimorphism.
		\item [(3)] If $H_{a+1,b+1}(\includegraphics[width=0.025\textheight]{"b_smoothing"})=0=H_{a-1,b+1}(\includegraphics[width=0.025\textheight]{"b_smoothing"})$ then $\beta_*: H_{a+1,b+1}(\includegraphics[width=0.025\textheight]{"crossing"}) \to H_{a,b}(\includegraphics[width=0.025\textheight]{"a_smoothing"})$ is an isomorphism.
	\end{enumerate}
\end{theorem}

\section{Khovanov Homology of Torus links $T(2,n)$}\label{KHtoruslinks}
Torus links of the type $T(2,n)$ have always enjoyed importance regarding research in knot theory. Indeed, the KH of this type of links was first computed by M. Khovanov \cite{Kho1}. Subsequently, J. H. Przytycki provided a different approach to this computation by observing the connection of Khovanov homology and Hochschild homology \cite{Prz}.

\

In this section, we present an alternative, deviceful method of calculating the homology of torus links $T(2,n)$ by explicitly using the long exact sequence of Khovanov homology.

\begin{definition}\label{Toruslinksdefi}
	A \textbf{torus link} of type (p,q), also referred to as a $(p,q)-$torus link is a link ambient isotopic to a curve contained in a standard torus $T^{2}$. This curve wraps $p$ times around the longitude and it wraps $q$ times around the meridian. If $p$ and $q$ are relatively prime, then the torus link has only one component and therefore it is called a \textbf{torus knot}.
\end{definition}

Let $v= \textbf{\KPB}$ be a fixed crossing. Observe that a vertical smoothing at a crossing $v$ of a torus link $T(2,n)$ produces the trivial knot with $1-n$ twists.\\

\begin{figure}[h]
	\centering
	\includegraphics[width=1\linewidth]{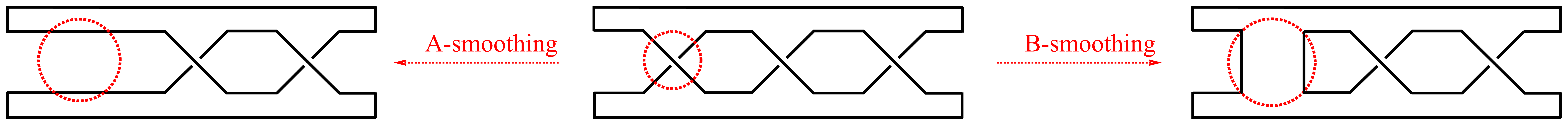}
	\caption{Torus knot $T(2,n)=T(2,3)$ (center), Torus knot $T(2,n-1)=T(2,2)$ (left), and the trivial knot with framing $1-n=-2$ (right).}
	\label{fig:torus-23}
\end{figure}

We now write the long exact sequence of KH. First, notice that maintaining the crossing (trivially) does not change the torus link $T(2,n)$, the horizontal smoothing of the crossing results in $T(2,n-1)$, and the vertical smoothing of the crossing yields the trivial knot with framing twisted $1-n$ times. In Figure \ref{fig:torus-23} this setting is illustrated for the torus knot $T(2,3)$.

The following theorem states the Khovanov homology of torus links $T(2,n)$ with $n>0$. The long exact sequence of KH plays a key role in the proof.
\begin{theorem} 
	Let $T(2,n)$ be a torus link with $n>0$. Its Khovanov homology $H_{a,b}(T(2,n))$ is given by:	
	\[ 
	H_{a,b}(T(2,n))=  \left\{   \begin{array}{ll}
		\mathbb{Z} & \mbox{for $(a,b)=(n,n)$ or $(a,b)=(-n,-3n)$},\\
		\mathbb{Z} & \mbox{for $a=n-2s$, $b=n-4s+4$ where $s$ is even and $0\leq s \leq n$,} \\
		\mathbb{Z} & \mbox{for $a=n-2s$. $b=n-4s$ where $s$ is odd and $3\leq s \leq n$,} \\
		\mathbb{Z}_2   & \mbox{for $a=n-2s$, $b=n-4s+4$ where $s$ odd and $3\leq s \leq n$,}\\
		0           & \mbox{otherwise.}
	\end{array}
	\right.
	\]

\end{theorem}

\begin{proof}
	The proof proceeds by induction on $n$.\\
	
	For the trivial knot $T(2,1)$, the theorem holds as only $H_{1,1}(T(2,1))= H_{1,5}(T(2,1))=\mathbb{Z}$ and the other homology groups are trivial.\\
	
	For the Hopf link $T(2,2)$ the theorem holds as can be verified with the Khovanov homology Table \ref{KHHopf}\footnote{A well-known computation. See, for instance, \cite{MV}.}.	
	
	\begin{table}[H] 
		\begin{minipage}{1\textwidth}
			\centering	
			\begin{tabular}{c||c|c|c|c|c} 
				\ \ $\textbf{b}  \ \ | \ \ \textbf{a}$ \ \  & \ \ \ -2 \ \ \ & \ \ \ 0 \ \ \ & \ \ \ 2 \ \ \  \\
				\hline \hline 
				6                                           &               &               &   $\mathbb{Z}$   \\ 
				\hline
				2                                           &               &                & $\mathbb{Z}$     \\
				\hline
				-2                                           &   $\mathbb{Z}$ &               &              \\
				\hline
				-6                                           & $\mathbb{Z}$  &                 &               \\
				\hline			
				
			\end{tabular}
		\end{minipage}
		\caption{Khovanov homology table of the torus knot $T(2,2)$.}  \label{KHHopf}
	\end{table}
	
	As the inductive hypothesis, suppose the result holds for $n-1$ where $n > 2$. Furthermore, consider the case where the map $\beta_{\ast}: H_{a,b}(T(2,n)) \longrightarrow H_{a-1,b-1}(T(2,n-1))$ is not necessarily an isomorphism. Thus, the long exact sequence should be closely analyzed. First, recall that after solving the crossing $v$ with a $B$-smoothing, the trivial knot with framing $1-n$ is obtained, i.e. $T_{B}(2,n)=\bigcirc^{1-n}$. Hence, its homology is given by: 
	
	\[ H_{x,y}(\bigcirc^{1-n})=
	\left\{ \begin{array}{ll}
		\mathbb{Z} & \mbox{when $(x,y)=(1-n,3(1-n)\pm 2)$},\\
		0  & \mbox{otherwise}.
	\end{array}
	\right.
	\]
	
	Therefore, the map $\beta_{\ast}$ is not necessarily an isomorphism if $H_{x,y}(\bigcirc^{1-n}) \neq 0$. More precisely, $H_{x,y}(\bigcirc^{1-n}) \neq 0$ when $(x,y)=(1-n,3(1-n)-2)=(1-n,1-3n)$ or $(x,y)=(1-n,3(1-n)+2))=(1-n,5-3n).$
	
\
	
	\hspace*{0.5cm}\textbf{Case (i)} Let $(x,y)=(1-n,1-3n)$:\\
	
	 In the long exact sequence, the neighborhood of $H_{1-n,1-3n}(\bigcirc^{1-n})$ has the form:
	\begin{equation*}
		\begin{split}
			\xrightarrow{\makebox[1cm]{}} H_{1-n,-1-3n}(T(2,n-1)) &   \xrightarrow{\makebox[2cm]{$  \textcolor{blue}{\left( \partial^{Conn} \right)}$}} H_{1-n,1-3n}(\bigcirc^{1-n}) \xrightarrow{\makebox[1.8cm]{\textcolor{blue}{$\alpha_{\ast}$}}}\\		 
			H_{-n,-3n}(T(2,n))\xrightarrow{\makebox[1.5cm]{\textcolor{blue}{$\beta_{\ast}$}}} & H_{-n-1,-3n-1}(T(2,n-1)) \xrightarrow{\makebox[1cm]{}}. \\
		\end{split}
	\end{equation*}	
	
	By the inductive hypothesis, $H_{1-n,-1-3n}(T(2,n-1))$ and $H_{-n-1,-3n-1}(T(2,n-1))$ are trivial. Hence, the previous long exact sequence becomes:
	\begin{equation*}
		\begin{split}
			\xrightarrow{\makebox[1cm]{}} 0 \xrightarrow{\makebox[2cm]{$  \textcolor{blue}{\left( \partial^{Conn} \right)}$}} H_{1-n,1-3n}(\bigcirc^{1-n}) & \xrightarrow{\makebox[1.8cm]{\textcolor{blue}{$\alpha_{\ast}$}}}\\		 
			H_{-n,-3n}(T(2,n))\xrightarrow{\makebox[1.5cm]{\textcolor{blue}{$\beta_{\ast}$}}} & 0 \xrightarrow{\makebox[1cm]{}}. \\
		\end{split}
	\end{equation*}		
	
	Then, $H_{1-n,1-3n}(\bigcirc^{1-n})=H_{-n,-3n}(T(2,n))=\mathbb{Z}$, and so the theorem holds as $(a,b)=(-n,-3n)$.\\
	
	\hspace*{0.5cm}\textbf{Case (ii)} Let $(x,y)=(1-n,5-3n)$:\\
	
	 In the long exact sequence, the neighborhood $H_{1-n,5-3n}(\bigcirc^{1-n})$ has the form:
	\begin{equation*}
		\begin{split}
			0 \xrightarrow{\makebox[1cm]{}} H_{2-n,4-3n}(T(2,n)) \xrightarrow{\makebox[1.5cm]{\textcolor{blue}{$\beta_{\ast}$}}} & H_{1-n,3-3n}(T(2,n-1))	 
			\xrightarrow{\makebox[2cm]{$  \textcolor{blue}{\left( \partial^{Conn} \right)}$}}\\ H_{1-n,5-3n}(\bigcirc^{1-n})\xrightarrow{\makebox[1.5cm]{\textcolor{blue}{$\alpha_{\ast}$}}} H_{-n,4-3n}(T(2,n))& \xrightarrow{\makebox[1cm]{\textcolor{blue}{$\beta_{\ast}$}}} \\
			&H_{-n-1,3-3n}(T(2,n-1)) \xrightarrow{\makebox[1cm]{}}. \\
		\end{split}
	\end{equation*}	
	
	By the inductive hypothesis, $H_{-n-1,3-3n}(T(2,n-1))=0$ and $H_{1-n,3-3n}(T(2,n-1))=H_{-(n-1),-3(n-1)}=\mathbb{Z}$. Thus, the previous long exact sequence becomes:
	\begin{equation*}
		\begin{split}
			0 \xrightarrow{\makebox[1cm]{}} H_{2-n,4-3n}(T(2,n)) \xrightarrow{\makebox[1.5cm]{\textcolor{blue}{$\beta_{\ast}$}}}  \mathbb{Z}	 
			\xrightarrow{\makebox[2cm]{$  \textcolor{blue}{\left( \partial^{Conn} \right)}$}} & \\ \mathbb{Z} \xrightarrow{\makebox[1.5cm]{\textcolor{blue}{$\alpha_{\ast}$}}} H_{-n,4-3n}(T(2,n))& \xrightarrow{\makebox[1cm]{\textcolor{blue}{$\beta_{\ast}$}}} 
			0 \xrightarrow{\makebox[1cm]{}}. \\
		\end{split}
	\end{equation*}	
	
	To determine the other entries of the long exact sequence we need to understand the map $\partial^{Conn}_{1-n,3-3n}: \mathbb{Z} \longrightarrow \mathbb{Z}$. \\
	There are two general possibilities: the map is the zero map or it is multiplication by $k>0$. \\
	
	First suppose the map $\partial^{Conn}_{1-n,3-3n}$ is the zero map. In this case,
	$$H_{2-n,4-3n}(T(2,n))=H_{1-n,3-3n}(T(2,n-1))=\mathbb{Z} \ \ \ \ \ and,$$
	$$H_{-n,4-3n}(T(2,n))=H_{1-n,5-3n}(\bigcirc^{1-n})=\mathbb{Z}.$$
	Then, $H_{2-n,4-3n}(T(2,n))=\mathbb{Z}$ and in this case we have that $(a,b)=(2-n,4-3n)=(n-2n+2,n-4n+4)=(n-2(n-1),n-4(n-1))=(n-2s,n-4s)$ for $s=n-1$. Similarly, $H_{-n,4-3n}(T(2,n))=\mathbb{Z}$ where $(a,b)=(-n,4-3n)=(n-2s,n-4s+4)$ for $s=n$ and thus, the theorem holds.\\
	
	Finally, suppose the map $\partial^{Conn}_{1-n,3-3n}$ is multiplication by $k>0$. This case was addressed in 2006 by  M. Pabiniak, J. H. Przytycki, and R. Sazdanović in \cite{PPS} by studying the algebra of truncated polynomials $\mathcal{A}_{m}=\quotient{\mathbb{Z}[x]}{(x^{m})}$, which for $m=2$ is closely related to Khovanov homology \cite{Kho1, BN}. More precisely, when $n$ is even, the map $\partial^{Conn}_{1-n,3-3n}$ is the zero map and thus we have  $H_{2-n,4-3n}(T(2,n))=\mathbb{Z}=H_{-n,4-3n}(T(2,n))$ as noted in the previous part. On the other hand, when $n$ is odd, the map $\partial^{Conn}_{1-n,3-3n}$ is multiplication by $k=2$ and we have that $H_{2-n,4-3n}(T(2,n))=0$, $H_{-n,4-3n}(T(2,n))=\mathbb{Z}_{2}$, and so the theorem holds.
	
\end{proof}

\section{Example: KH tables for T(2,11) and T(2,12)}\label{examplesKHtables}
We present in this section the KH Table \ref{Table 1} and KH Table \ref{Table 2} for the torus links $T(2,11)$ and $T(2,12)$, respectively.

\

Observe that $H_{a,b}(T(2,n))$ has support on two diagonals containing $(n, n)$ or $(n,n+4)$. In other words, $H_{a,b}(T(2,n))$ is nontrivial only for $(a,b)=(n-2s,n-4s)$ or $(a,b)=(n-2s,n-4s +4)$. Moreover, $H_{a,b}(T(2,n))$ has torsion only in groups for which $(a,b)=(n-2s,n-4s +4)$.

\begin{table}[ht] 
	\begin{minipage}{1\textwidth}
		\centering	
		\begin{tabular}{c||c|c|c|c|c|c|c|c|c|c|c|c|}
			\ \ $\textbf{b} \ \ | \  \ \textbf{a}$  &-11&-9& -7& -5  &  -3  &  -1  &  1 & 3 & 5 & 7& 9 & 11 \\
			\hline \hline
			15                 & & & & & & & & & & & & $\mathbb{Z}$  \\
			\hline
			11                 & & & & & & & & & & & & $\mathbb{Z}$    \\
			\hline
			7                  & & & & & & & & & &     $\mathbb{Z}$ & &  \\
			\hline
			3                  & & & & & & & & &       $\mathbb{Z}_{2}$ & & &  \\
			\hline
			-1                  & & & & & & & &            $\mathbb{Z}$  &  $\mathbb{Z}$ & & & \\
			\hline
			-5                 & & & & & & &  $\mathbb{Z}_{2}$ & & & & & \\
			\hline
			-9                 & & & & & &  $\mathbb{Z}$ & $\mathbb{Z}$ & & & & &    \\
			\hline
			-13                & & & & &   $\mathbb{Z}_{2}$ & & & & & & & \\
			\hline
			-17                & & & &  $\mathbb{Z}$ & $\mathbb{Z}$ & & & & & & &    \\
			\hline
			-21                & & &   $\mathbb{Z}_{2}$ & & & & & & & & & \\
			\hline
			-25                & &  $\mathbb{Z}$ & $\mathbb{Z}$ & & & & & & & & &    \\
			\hline
			-29                &    $\mathbb{Z}_{2}$ & & & & & & & & & & & \\
			\hline
			-33               & $\mathbb{Z}$   &&&&&&&&&&& \\
			\hline
		\end{tabular}
	\end{minipage}
	\caption{Khovanov homology of the torus knot $T(2,11)$.} \label{Table 1}
\end{table}

\begin{table}[ht]
	\begin{minipage}{1\textwidth}
		\centering
		\begin{tabular}{c||c|c|c|c|c|c|c|c|c|c|c|c|c|}
			\ \ $\textbf{b} \ \ | \  \ \textbf{a}$ &-12 &-10&-8& -6& -4  &  -2  &  0  &  2 & 4 & 6 & 8& 10 & 12 \\
			\hline \hline
			16               & & & & & & & & & & & & & $\mathbb{Z}$  \\
			\hline
			12               & & & & & & & & & & & & & $\mathbb{Z}$    \\
			\hline
			8                & & & & & & & & & & &     $\mathbb{Z}$ & &  \\
			\hline
			4                & & & & & & & & & &       $\mathbb{Z}_{2}$ & & &  \\
			\hline
			0               &  & & & & & & & &            $\mathbb{Z}$  &  $\mathbb{Z}$ & & & \\
			\hline
			-4               & & & & & & & &  $\mathbb{Z}_{2}$ & & & & & \\
			\hline
			-8               & & & & & & &  $\mathbb{Z}$ & $\mathbb{Z}$ & & & & &    \\
			\hline
			-12              & & & & & &   $\mathbb{Z}_{2}$ & & & & & & & \\
			\hline
			-16              & & & & &  $\mathbb{Z}$ & $\mathbb{Z}$ & & & & & & &    \\
			\hline
			-20              & & & &   $\mathbb{Z}_{2}$ & & & & & & & & & \\
			\hline
			-24              & & &  $\mathbb{Z}$ & $\mathbb{Z}$ & & & & & & & & &    \\
			\hline
			-28              & &    $\mathbb{Z}_2$ & & & & & & & & & & & \\
			\hline
			-32               & $\mathbb{Z}$   & $\mathbb{Z}$ &&&&&&&&&&& \\
			\hline
			-36              &  $\mathbb{Z}$   &&&&&&&&&&&& \\
			\hline
		\end{tabular}
	\end{minipage}
	\caption{Khovanov homology of the torus link $T(2,12)$.}\label{Table 2}
\end{table}	

\section{Future directions}\label{Futuredirec} 

Khovanov Homology has been computed for many links and the experimental evidence shows abundance of $2-$torsion while other torsion groups appear (until recently) very seldom; see for instance \cite{Kho1, Shu}. In \cite{AP} the authors studied the notion of adequate diagrams. They observed that for an $n-$crossing $A-$adequate link diagram $D$, the highest grading in KH is given by $H_{n,\ast}(D)=H_{n, n+2|s_{A}|}(D)=\mathbb{Z}$, where $S_{A}$ denotes the Kauffman state with an $A$ marker assigned at all its crossings. Moreover, they showed that the next nontrivial homology group $H_{n-2,n+2|s_{A}|-4}(D)$, has $\mathbb{Z}_{2}-$torsion as long as the graph associated to the state $s_{A}$, $G(D)=G_{s_{A}}(D)$, is not bipartite. Furthermore, the group $H_{n-2,n+2|s_{A}|-4}(D)$ was explicitly computed for $A-$adequate links in \cite{PPS}. In particular, for a connected $A-$adequate diagram $D$,
$$tor\left( H_{n-2,n+2|s_{A}|-4}(D) \right) = \left\lbrace 
\begin{array}{ll}
	\mathbb{Z}_{2}, & for \ \ G(D) \ \ having \ \ and \ \ odd \ \ cycle \\
	0, & for \ \ a \ \ bipartite \ \ graph.
\end{array}
\right.$$
They additionally showed that for a strongly $A-$adequate diagram $D$ with the graph $G(D)$ containing an even cycle, this group contains $\mathbb{Z}_{2}-$torsion. An idea for further research might be to analyze whether the long exact sequence of KH can be used in full generality to compute the KH for this type of links and to detect torsion, extending the work in \cite{PS}, for instance. In the fairly recent papers \cite{SS1} and \cite{SS2}, the authors further extended on the results from \cite{PPS} and \cite{PS}; for instance, in the latest paper they show the triviality of Khovanov homology in certain gradings and study the extreme Khovanov groups and their corresponding coefficients of the Jones polynomial. These papers elaborate on the close relation between Khovanov homology and chromatic graph homology. 

\

Lately there have been significant developments in the study of torsion in Khovanov homology, see for instance \cite{PBIMW, Muk, MS}. The method used in the calculation of the KH of torus links $T(2,n)$ can be used to construct links having torsion different from $\mathbb{Z}_{2}$ in its Khovanov homology. It is a reasonable and promising idea to test the full power of the long exact sequence of KH and expand on the results from \cite{MPSWY}. For example, a natural research path is to adapt the long exact sequence to calculate the KH of other families of links. One can consider torus links of type $T(3,k)$ since a horizontal smoothing would produce a $T(2,n)$ torus link. Likewise, one could consider pretzel links with three columns, as an horizontal smoothing would produce a similar result. 

\section*{Acknowledgments}
The author proudly acknowledges the support of the National Science Foundation through Grant DMS-2212736. Furthermore, the author gives special thanks to his doctoral advisor, J\'ozef H. Przytycki, whose guidance has been essential in exploring the beautiful knot theory world.


\end{document}